
\input amstex.tex
\documentstyle{amsppt}
\magnification=\magstep1
\hsize=12.5cm
\vsize=18cm
\hoffset=1cm
\voffset=2cm
\def\DJ{\leavevmode\setbox0=\hbox{D}\kern0pt\rlap
 {\kern.04em\raise.188\ht0\hbox{-}}D}
\footline={\hss{\vbox to 2cm{\vfil\hbox{\rm\folio}}}\hss}
\nopagenumbers 
\font\ff=cmr8
\def\txt#1{{\textstyle{#1}}}
\baselineskip=13pt
\def\hf{{\textstyle{1\over2}}}
\def\a{\alpha}
\def\d{{\,\roman d}}
\def\e{\varepsilon}

\def\G{\Gamma}
\def\k{\kappa}
\def\s{\sigma}

\def\={\;=\;}
\def\zx{\zeta(\hf+ix)}
\def\zt{\zeta(\hf+it)}

  \def\kk{\kappa_h}
\def\R{\Re{\roman e}\,} \def\I{\Im{\roman m}\,} \def\s{\sigma}
\def\z{\zeta}

\def\H{H_j^3({\txt{1\over2}})} 
\font\teneufm=eufm10
\font\seveneufm=eufm7
\font\fiveeufm=eufm5
\newfam\eufmfam
\textfont\eufmfam=\teneufm
\scriptfont\eufmfam=\seveneufm
\scriptscriptfont\eufmfam=\fiveeufm
\def\mathfrak#1{{\fam\eufmfam\relax#1}}

\font\tenmsb=msbm10
\font\sevenmsb=msbm7
\font\fivemsb=msbm5
\newfam\msbfam
\textfont\msbfam=\tenmsb
\scriptfont\msbfam=\sevenmsb
\scriptscriptfont\msbfam=\fivemsb
\def\Bbb#1{{\fam\msbfam #1}}

\def \NN {\Bbb N}
\def \CC {\Bbb C}

\def \ZZ {\Bbb Z}

\def\rightheadline{{\hfil{\ff
Laplace transform of the fourth moment of $|\zx|$}\hfil\tenrm\folio}}

\def\leftheadline{{\tenrm\folio\hfil{\ff
Aleksandar Ivi\'c }\hfil}}
\def\emptyheadline{\hfil}
\headline{\ifnum\pageno=1 \emptyheadline\else
\ifodd\pageno \rightheadline \else \leftheadline\fi\fi}

\topmatter
\title THE LAPLACE TRANSFORM OF THE FOURTH MOMENT
OF THE ZETA-FUNCTION \endtitle
\author   Aleksandar Ivi\'c 
\bigskip
\endauthor
\dedicatory
Univ. Beograd. Publ. Elektrotehn. Fak. Ser. Mat. 
{\bf11}${\roman(2000), 41-48}$.
\enddedicatory
\address{
Aleksandar Ivi\'c, Katedra Matematike RGF-a
Universiteta u Beogradu, \DJ u\v sina 7, 11000 Beograd,
Serbia (Yugoslavia). 
}
\endaddress
\email {\tt aleks\@ivic.matf.bg.ac.yu, eivica\@ubbg.etf.bg.ac.yu} \endemail
\abstract
{The Laplace transform of $|\zx|^4$ is investigated, for which a precise
expression is obtained, valid in a certain region in the
complex plane. The method of proof is based on complex integration and
spectral theory of the non-Euclidean Laplacian.}
\endabstract
\endtopmatter

\heading {\bf 1. INTRODUCTION} \endheading

Laplace transforms play an important r\^ole in analytic number
theory. Of special interest in the theory of the Riemann
zeta-function $\z(s)$ are the Laplace transforms
$$
L_k(s) \;:=\; \int_0^\infty |\zx|^{2k}e^{-sx}\d x\qquad(k \in \NN,\,\R s > 0).
\leqno(1.1)
$$
E.C. Titchmarsh's well-known monograph [20, Chapter 7] gives a discussion
of $L_k(s)$ when $s = \s$ is real and $\s \to 0+$, especially detailed
in the cases $k=1$ and $k=2$. Indeed, a classical result of H. Kober
[14] says that, as $\s \to 0+$,
$$
L_1(2\s) = {\gamma-\log(4\pi\s)\over2\sin\s} +
 \sum_{n=0}^Nc_n\s^n + O(\s^{N+1})
\leqno(1.2)
$$
for any given integer $N \ge 1$, where the $c_n$'s are effectively
computable constants and $\gamma = 0.577\ldots\,$ is Euler's constant.
For complex values of $s$ the function $L_1(s)$ was studied by
F.V. Atkinson [1], and more recently by M. Jutila [13], who noted that
Atkinson's argument gives
$$
L_1(s) = -ie^{{1\over2}is}(\log(2\pi)-\gamma + ({\pi\over2}-s)i)
+ 2\pi e^{-{1\over2}is}\sum_{n=1}^\infty
d(n)\exp(-2\pi ine^{-is}) + \lambda_1(s) \leqno(1.3)
$$
in the strip $0 < \R s < \pi$, where the function $\lambda_1(s)$ is 
holomorphic in the strip $|\R s| < \pi$. Moreover, in any strip
$|\R s| \le \theta$ with $0 < \theta < \pi$, we have
$$
\lambda_1(s) \;\ll_\theta \;(|s|+1)^{-1}.
$$
In [12] M. Jutila gave a discussion on the application of Laplace
transforms to the evaluation of sums of coefficients of certain
Dirichlet series.

F.V. Atkinson [2] obtained the asymptotic formula
$$
L_2(\s) = {1\over\s}(A\log^4{1\over\s} + B\log^3{1\over\s}
+ C\log^2{1\over\s} + D\log {1\over\s} + E) + \lambda_2(\s),\leqno(1.4)
$$
where $\s \to 0+$,
$$
A = {1\over2\pi^2},\,B =\pi^{-2}(2\log(2\pi) - 6\gamma + 24\z'(2)\pi^{-2})
$$
and
$$
 \lambda_2(\s) \;\ll_\e\;\left({1\over\s}\right)^{{13\over14}+\e}.\leqno(1.5)
$$
He also indicated how, by the use of estimates for Kloosterman sums,
one can improve the exponent ${13\over14}$ in (1.5) to ${8\over9}$. This
is of historical interest, since it is one of the first instances of
application of Kloosterman sums to analytic number theory. Atkinson
in fact showed that ($\s = \R s > 0$)
$$
L_2(s) \;=\;4\pi e^{-{1\over2}s}\sum_{n=1}^\infty d_4(n)
K_0(4\pi i\sqrt{n}e^{-{1\over2}s}) + \phi(s),\leqno(1.6)
$$
where $d_4(n)$ is the divisor function generated by $\z^4(s)$, $K_0$
is the Bessel function, and the series in (1.6) as well as $\phi(s)$
are both analytic in the region $|s| < \pi$. When $s = \s \to 0+$ one
can use the asymptotic formula
$$
K_0(z) \= \hf\sqrt{\pi}z^{-1/2}e^{-z}\left(1 - 8z^{-1} + O(|z|^{-2})\right)
\quad(|\arg z| < \theta < {3\pi\over2},\,|z| \ge 1)
$$
and then, by delicate analysis, one can deduce (1.4)--(1.5) from (1.6).

The author [5] gave explicit, albeit complicated expressions for
the remaining coefficients $C,D$ and $E$ in (1.4). More importantly,
he applied a result on the fourth moment of $|\zt|$, obtained jointly
with Y. Motohashi [9], [11] (see also [4]), to establish that
$$
\lambda_2(\s) \;\ll\; \s^{-1/2}\qquad(\s\to 0+),\leqno(1.7)
$$
and this is where the matter presently rests.

\medskip
For $k \ge 3$ not much is known about $L_k(s)$, even when $s = \s \to 0+$.
This is not surprising, since not much is known about upper bounds for
$$
I_k(T) \;:=\; \int_0^T|\zt|^{2k}\d t\qquad(k \ge 3,\, k \in \NN).
$$
For a discsussion on $I_k(T)$ the reader is referred to the author's
monographs [3] and [4]. One trivially has
$$
I_k(T) \le e\int_0^\infty|\zt|^{2k}e^{-t/T}\d t = eL_k({1\over T}).
\leqno(1.8)
$$
Thus any nontrivial bound of the form
$$
L_k(\s) \;\ll_\e\; \left({1\over\s}\right)^{c_k+\e}\qquad(\s\to 0+,\,
c_k \ge 1)\leqno(1.9)
$$
gives, in view of (1.8) ($\s = 1/T$), the bound
$$
I_k(T) \;\ll_\e\; T^{c_k+\e}.\leqno(1.10)
$$
Conversely, if (1.10) holds, then we obtain (1.9) from the identity
$$
L_k({1\over T}) \= {1\over T}\,\int_0^\infty I_k(t)e^{-t/T}\d t,
$$
which is easily established by integration by parts.

\bigskip\medskip
\heading {\bf  2. SPECTRAL THEORY AND THE 
LAPLACE TRANSFORM of $|\zx|^4$} \endheading

\bigskip\medskip
 
The purpose of this paper is to consider $L_2(s)$, where $s$ is a complex
variable, and to prove a result analogous to (1.3), valid in a certain
region in $\CC$. We shall not use Atkinson's method and try to elaborate
on (1.6).  Our main tools are  powerful methods from spectral
theory, by which recently much advance has been made in connection
with $I_2(T)$. For a competent and extensive account of
spectral theory the reader is referred to Y. Motohashi's monograph [19].
Some of the relevant papers on $I_2(T)$ are [6]--[11], [14]--[18] and [21].

We begin by stating briefly  the necessary notation 
involving the spectral theory of the non-Euclidean Laplacian.
As usual $\,\{\lambda_j = \kappa_j^2 + {1\over4}\} \,\cup\, \{0\}\,$ will
denote the discrete spectrum of the non-Euclidean Laplacian acting
on $\,SL(2,\ZZ)\,$--automorphic forms, and
$\a_j = |\rho_j(1)|^2(\cosh\pi\kappa_j)^{-1}$, where
$\rho_j(1)$ is the first Fourier coefficient of the Maass wave form
corresponding to the eigenvalue $\lambda_j$ to which the Hecke series
$H_j(s)$ is attached. We note that
$$
\sum_{\k_j\le K}\a_jH_j^3(\hf) \ll K^2\log^CK \qquad(C > 0).\leqno(2.1)
$$
Our result is the following

\medskip
THEOREM. 
{\it Let $0 \le \phi < {\pi\over2}$ be given. Then for $0 < |s| \le 1$ and
$|\arg s| \le \phi$ we have}
$$\eqalign{&
L_2(s) = {1\over s}(A\log^4{1\over s} + B\log^3{1\over s} +
C\log^2{1\over s} + D\log{1\over s} + E) + G_2(s)\cr&
+ s^{-{1\over2}}\left\{\sum_{j=1}^\infty \a_j\H\left(
s^{-i\k_j}R(\k_j)\G(\hf + i\k_j) + 
s^{i\k_j}R(-\k_j)\G(\hf - i\k_j) \right)\right\}
,\cr}\leqno(2.2)
$$
{\it where}
$$
R(y) \;:=\; \sqrt{{\pi\over2}}{\Bigl(2^{-iy}{\G({1\over4} - {i\over2}y)
\over\G({1\over4} + {i\over2}y)}\Bigr)}^3\G(2iy)\cosh(\pi y)
\leqno(2.3)
$$
{\it and in the above region $G_2(s)$ is a regular function satisfying
($C > 0$ is a suitable constant)}
$$
G_2(s) \ll |s|^{-1/2}\exp\left\{
-{C\log(|s|^{-1}+20)\over(\log\log(|s|^{-1}+20))^{2/3}
(\log\log\log(|s|^{-1}+20))^{1/3}}\right\}.
\leqno(2.4)
$$

{\bf Remark 1}. The constants $A,\,B,\,C,\,D,\,E$ in (2.2) are the same
ones as in (1.4).

{\bf Remark 2}. From Stirling's formula for the gamma-function it
follows that $R(\k_j) \ll \k_j^{-1/2}$. In view of (2.1) 
this means that the series in (2.2) is absolutely convergent and
uniformly bounded in $s$ when $s = \s$ is real. Therefore, when
$s = \s \to 0+$, (2.2) gives a refinement of (1.7).

{\bf Remark 3}. From (1.4) and (1.7) it transpires that 
$\lambda(\s)$ is an error term when $0 < \s < 1$. For this reason we
considered the values $0 < |s| \le 1$ in (2.2), although one could treat
the case $|s| > 1$ as well.

{\bf Remark 4}. From (2.2) and elementary properties of the 
Laplace transform one can easily obtain the Laplace transform
of
$$
E_2(T) := \int_0^T|\zt|^4\d t - TP_4(\log T),\quad
P_4(x) = \sum_{j=0}^4a_jx^j,
$$
where $a_4 = 1/(2\pi^2)$ (for the evaluation of the remaining $a_j$'s,
see [5]).
\medskip

\heading{\bf 3. PROOF OF THE THEOREM}\endheading

Note first that
in the integral defining $L_2(s)$ it suffices consider only the range
of integration $\,[1,\,\infty\,]\,$, since the range $\,[0,\,1]\,$
trivially contributes $\ll 1$.
We start from the well-known integral
$$
e^{-z} = {1\over2\pi i}\int_{(c)}\G(s)z^{-s}\d s\qquad(\R z > 0, c > 0),
\leqno(3.1)
$$
and the function
$$
{\Cal Z}_2(w) := \int_1^\infty|\zt|^4t^{-w}\d t\quad(\R w > 1),
$$
introduced and studied by Y. Motohashi [17], [19]. Here as usual
$$
\int_{(\s)} \;=\;\lim_{T\to\infty}\int_{\s-iT}^{\s+iT}.
$$ 
Y. Motohashi
showed that ${\Cal Z}_2(s)$ has meromorphic continuation over $\CC$.
In the half-plane $\s = \R s > 0$
it has the following singularities: the pole $s = 1$ of order 5,
simple poles at $s = {1\over2} \pm i\k_j\,(\k_j
= \sqrt{\lambda_j - {1\over4}})$ and poles at $s = \hf\rho$, 
where $\rho$ denotes complex
zeros of $\z(s)$. The residue of ${\Cal Z}_2(s)$ at 
$s = {1\over2} + i\kk$ equals
$$
R_0(\kk) = \sqrt{{\pi\over2}}{\Bigl(2^{-i\kk}{\G({1\over4} - {i\over2}\kk)
\over\G({1\over4} + {i\over2}\kk)}\Bigr)}^3\G(2i\kk)\cosh(\pi\kk)
\sum_{\k_j=\kk}\alpha_j \H, 
$$
and the residue at $s = {1\over2} - i\kk$ equals $\overline{R_0(\kk)}$.

From (3.1) we have, for $c > 1$, $0 < |s| \le 1$ and
$|\arg s| \le \phi$,
$$
\eqalign{&
\int_1^\infty|\zx|^4e^{-sx}\d x\cr&
= \int_1^\infty|\zx|^4\left(
{1\over2\pi i}\int_{(c)}\G(w)(sx)^{-w}\d w\right)\d x\cr&
= {1\over2\pi i}\int_{(c)}\G(w)s^{-w}{\Cal Z}_2(w)\d w.\cr}\leqno(3.2)
$$
In (3.2) we shift the line of integration to the contour ${\Cal L}$
($w = u + iv$) consisting of the curves
$$
u = \hf - C\log^{-2/3}|v|(\log\log|v|)^{-1/3}\qquad(|v| \ge v_0 > 0,\,
C > 0)\leqno(3.3)
$$
and the segment
$$
u = u_0,\,u_0 = \hf - C\log^{-2/3}|v_0|(\log\log|v_0|)^{-1/3}, \;|v| \le v_0.
$$
Namely $\z(s) \not= 0$ (see e.g., [3, Chapter 6]) for
$$
\s \ge 1 - A(\log t)^{-2/3}(\log\log t)^{-1/3}
\qquad(s = \s + it,\, t \ge t_0 > 0,\,A > 0)
$$
and a suitable constant $A$.
The function ${\Cal Z}_2(w)$ will be regular on ${\Cal L}$,
since the poles $w = \hf\rho,\,\z(\rho) = 0)$ lie to the left of ${\Cal L}$.
For any given $\eta > 0$ one  has
$$
{\Cal Z}_2(w) \;\ll\;e^{\eta|\I w|}\qquad(w \in {\Cal L})\leqno(3.4)
$$
if $C$ in (3.3) is taken sufficiently small. This follows easily from
the proof of Lemma 1 of [7], which shows that the order of ${\Cal Z}_2(w)$
is of the order given by (3.4) if $\R w > 0$ and $w$ stays away from the
poles of ${\Cal Z}_2(w)$. In the forthcoming work [8] it is even shown
that in the above region ${\Cal Z}_2(w)$ is of polynomial growth
in $|v|$, which is more than what is required for our present purpose.
Thus by the residue theorem we obtain from (3.2)
$$
\int_1^\infty|\zx|^4e^{-sx}\d x = \sum\,
{\roman {Res}}\;\G(w)s^{-w}{\Cal Z}_2(w)
+ {1\over2\pi i}\int_{\Cal L}\G(w)s^{-w}{\Cal Z}_2(w)\d w,\leqno(3.5)
$$
and the last integral is regular for $0 < |s| \le 1$ and
$|\arg s| \le \phi$.
There are residues at $w = 1$ and at $w = \hf \pm i\k_j$. The contribution
from $w = 1$ is (see (1.4))
$$
{1\over s}(A\log^4{1\over s} + B\log^3{1\over s} +
C\log^2{1\over s} + D\log{1\over s} + E),
$$
while the residues at $w = \hf \pm i\k_j$ yield
$$
s^{-1/2}\left\{\sum_{j=1}^\infty \a_j\H\left(
s^{-i\k_j}R(\k_j)\G(\hf+i\k_j) + 
s^{i\k_j}R(-\k_j)\G(\hf-i\k_j) \right)\right\},
$$
where $R(y)$ is defined by (2.3). 

Write
$$
\int_{\Cal L}\G(w)s^{-w}{\Cal Z}_2(w)\d w \= I_1 + I_2,
$$
say, where in $I_1$ we have $|v| \le V$, and in $I_2$ we have
$|v| > V$, where $V (\gg 1)$ is a parameter to be chosen. Since
$$
|s^{-w}| \= |s|^{-u}e^{|v\arg s|} \quad(|\arg s| \le \phi),
\;\,\G(w) \ll |v|^{u-{1\over2}}
e^{-{\pi\over2}|v|}\quad(v \gg 1),
$$
then setting $\delta(x) = C\log^{-2/3}x(\log\log x)^{1/3}$ we obtain
$$
I_1 \ll
|s|^{-u_0} + \int_{v_0}^V |s|^{-({1\over2}-\delta(v))}
e^{-({\pi\over2}-\phi)v}\d v
\ll |s|^{-u_0} + |s|^{-({1\over2}-\delta(V))}.
$$
Similarly we have
$$
I_2 \ll |s|^{-{1\over2}}\int_V^\infty e^{-({\pi\over2}- \phi)v}\d v
= |s|^{-{1\over2}}\,{ e^{-({\pi\over2}- \phi)V}\over ({\pi\over2}- \phi)}.
$$
Finally choosing
$$
V \= C_1\log(|s|^{-1}+20)
$$
with suitable $C_1 > 0$ we obtain ($C_2 > 0$)
$$\eqalign{&
\int_{\Cal L}\G(w)s^{-w}{\Cal Z}_2(w)\d w
\cr&\ll |s|^{-1/2}\exp\left\{
-{C_2\log(|s|^{-1}+20)\over(\log\log(|s|^{-1}+20))^{2/3}
(\log\log\log(|s|^{-1}+20))^{1/3}}\right\},\cr}
$$
which in view of (3.5)  completes the proof of the Theorem.

\bigskip\bigskip\bigskip
\Refs

\item{[1]} F.V. Atkinson, The mean value of the zeta-function on
the critical line, {\it Quart. J. Math. Oxford} {\bf 10}(1939), 122-128.

\item {[2]} F.V. Atkinson, The mean value of the zeta-function on
the critical line, {\it Proc. London Math. Soc.} {\bf 47}(1941), 174-200.

\item {[3]} A. Ivi\'c,  The Riemann zeta-function, {\it John Wiley
and Sons}, New York, 1985.

\item {[4]} A. Ivi\'c,  Mean values of the Riemann zeta-function,
LN's {\bf 82}, {\it Tata Institute of Fundamental Research}, 
Bombay, 1991 (distr. by Springer Verlag, Berlin etc.).

\item{ [5]}  A. Ivi\'c,   On the fourth moment of the Riemann
zeta-function, {\it Publs. Inst. Math. (Belgrade)}
{ \bf57(71)}(1995), 101-110.

\item{[6]} A. Ivi\'c,  The Mellin transform and the Riemann zeta-function, in
``{\it Proceedings of the Conference on Elementary and Analytic Number Theory 
(Vienna, July} 18-20, 1996)", Universit\"at Wien \& Universit\"at f\"ur
Bodenkultur, Eds. W.G. Nowak and J. Schoi{\ss}engeier, Vienna 1996, 112-127.

\item{[7]} A. Ivi\'c, On the error term for the fourth moment of the
Riemann zeta-function, {\it J. London Math. Soc.} 
{\bf60}(2)(1999), 21-32.

\item{[8]} A. Ivi\'c, M. Jutila and Y. Motohashi, The Mellin transform of
power moments of the zeta-function, {\it Acta Arithmetica} 
 {\bf95}(2000), 305-342.

\item{ [9]} A. Ivi\'c and Y. Motohashi,  A note on the mean value of
the zeta and L-functions VII, {\it  Proc. Japan Acad. Ser. A }
{\bf 66}(1990), 150-152.

\item{ [10]} A. Ivi\'c and Y. Motohashi,  The mean square of the
error term for the fourth moment of the zeta-function, {\it Proc. London Math.
Soc.} (3){\bf 66}(1994), 309-329.

\item {[11]} A. Ivi\'c and Y. Motohashi,  The fourth moment of the
Riemann zeta-function, {\it J. Number Theory} {\bf 51}(1995), 16-45.

\item{[12]}  M. Jutila, Mean values of Dirichlet series via Laplace 
transforms,
in {\it ``Analytic Number Theory"} (ed. Y. Motohashi), London Math. Soc.
 LNS {\bf247}, {\it Cambridge University Press}, Cambridge, 1997, 169-207.

\item{[13]}  M. Jutila, Atkinson's formula revisited, 
in {\it ``Voronoi's Impact on Modern Science"},
Book 1, Institute of  Mathematics, 
National Academy of Sciences of Ukraine, Kyiv, 1998, 139-154.

\item {[14]} H. Kober, Eine Mittelwertformel der Riemannschen Zetafunktion,
{\it Compositio Math.} {\bf3}(1936), 174-189.

\item {[15]} Y. Motohashi,  The fourth power mean of the Riemann
zeta-function, in {\it ``Proceedings of the Amalfi Conference on Analytic
Number Theory 1989"}, eds. E. Bombieri et al., Universit\`a di Salerno,
Salerno, 1992, 325-344.

\item{ [16]} Y. Motohashi,   An explicit formula for the fourth power 
mean of the Riemann zeta-function, {\it Acta Math. }{\bf 170}(1993), 181-220.

\item {[17]} Y. Motohashi,  A relation between the Riemann zeta-function
and the hyperbolic Laplacian, {\it Annali Scuola Norm. Sup. Pisa, Cl. Sci. IV
ser.} {\bf 22}(1995), 299-313.

\item {[18]} Y. Motohashi,  The Riemann zeta-function and the
non-Euclidean Laplacian, {\it Sugaku Expositions}, AMS {\bf 8}(1995), 59-87.

\item {[19]} Y. Motohashi,  Spectral theory of the Riemann
zeta-function, {\it Cambridge University Press}, Cambridge, 1997.

\item{[20]} E.C. Titchmarsh, The theory of the Riemann zeta-function
(2nd ed.), Clarendon Press, Oxford, 1986.

\item{ [21]}  N.I. Zavorotnyi, On the fourth moment 
of the Riemann zeta-function
(in Russian), in {\it ``Automorphic functions and number theory I"}, 
Coll. Sci. Works, Vladivostok, 1989, 69-125.

\bigskip
\bigskip
\bigskip
\bigskip
\bigskip
\bigskip

Aleksandar Ivi\'c

Katedra Matematike RGF-a

Universitet u Beogradu, \DJ u\v sina 7

11000 Beograd, Serbia (Yugoslavia)

\tt ivic\@rgf.bg.ac.yu, \enskip aivic\@matf.bg.ac.yu

\endRefs

\bye